\documentclass{article}
\usepackage{amsmath}
\usepackage{amsmath,amssymb,latexsym,color}
\usepackage[mathscr]{eucal}
\newtheorem{thm}{Theorem}[section]

\newtheorem{lemma}[thm]{Lemma}

\newtheorem{defin}[thm]{Definition}

\newtheorem{example}{Example}[section]

\newtheorem{remark}{Remark}[section]

\newcommand{\proof}{{\it Proof.\quad}}
\newcommand{\qed}{\hfill\Box\medskip}

\usepackage{CJK}
\begin{document}

\renewcommand{\baselinestretch}{1.2}

\title{\bf Erd\H{o}s-Ko-Rado theorem for vector spaces
 over residue class rings}

\author{
Jun Guo\thanks{Corresponding author. guojun$_-$lf@163.com}\\
{\footnotesize Department of Mathematics, Langfang Normal University, Langfang  065000,  China} }
 \date{ }
 \maketitle

\begin{abstract}
Let $h=\prod_{i=1}^{t}p_i^{s_i}$ be
its decomposition into a product of powers of distinct primes, and $\mathbb{Z}_{h}$ be the residue class ring modulo $h$.
Let  $\mathbb{Z}_{h}^{n}$ be the $n$-dimensional row vector space over $\mathbb{Z}_{h}$.
 A generalized Grassmann graph over $\mathbb{Z}_{h}$, denoted by $G_r(m,n,\mathbb{Z}_{h})$ ($G_r$ for short), has all
$m$-subspaces of $\mathbb{Z}_{h}^n$ as its vertices, and two distinct vertices
are adjacent if their intersection is of dimension $>m-r$,
where  $2\leq r\leq m+1\leq n$. In this paper, we determine the clique number
and geometric structures of maximum cliques of $G_r$. As a result, we obtain the Erd\H{o}s-Ko-Rado theorem for $\mathbb{Z}_{h}^{n}$.



\end{abstract}

\section{Introduction}
Let $\mathbb{Z}$ denote  the integer ring.
For $a,b,h\in \mathbb{Z}$, integers $a$ and $b$ are said to be {\it congruent modulo} $h$
 if $h$ divides $a-b$, and  denoted by $a\equiv b \mod h$.
 Suppose
  \begin{equation}\label{equa1}
  h=p_1^{s_1}p_2^{s_2}\cdots p_t^{s_t},
   \end{equation}
   where $p_i\,(1\leq i\leq t)$ are pairwise distinct primes, and $s_i\,(1\leq i\leq t)$ are positive integers.
  Let $\mathbb{Z}_{h}$ denote the residue class ring modulo $h$ and $\mathbb{Z}_{h}^\ast$ its unit group.
 Then $\mathbb{Z}_{h}$ is a  principal ideal ring and $|\mathbb{Z}_{h}^\ast|=h\prod_{i=1}^t(1-p_i^{-1}).$
  Let $[t]=\{1,2,\ldots,t\}$ and $J_{(\alpha_1,\alpha_2,\ldots,\alpha_t)}=(\prod_{i=1}^{t}p_i^{\alpha_i})$, where $0\leq \alpha_i\leq s_i$ for each $i\in[t]$. For brevity, we write $J_{(\alpha_1)}$ as $J_{\alpha_1}$ if $t=1$. By Lemma~2.3 in \cite{Guo3}, the principal ideals $J_{(\alpha_1,\alpha_2,\ldots,\alpha_t)}$, where $0\leq \alpha_{i}\leq s_i$ for each $i\in[t]$, are all the ideals of $\mathbb{Z}_{h}$.
For $a\in\mathbb{Z}$, we denote also by $a$ the congruence class of $a$ modulo $h$.

 Let  $\pi_i$ be the natural surjective homomorphism from
$(\mathbb{Z}_{h},+,\cdot)$  to $(\mathbb{Z}_{p_i^{s_i}},+,\cdot)$
for each $i\in[t]$.
Then there is an  isomorphism  as follows (cf. \cite{Ireland}):
\begin{equation}\label{equa2}
\begin{array}{rcl}
\pi:(\mathbb{Z}_{h},+,\cdot) &\rightarrow & (\mathbb{Z}_{p_1^{s_1}}\oplus \mathbb{Z}_{p_2^{s_2}}
 \oplus\cdots\oplus\mathbb{Z}_{p_t^{s_t}},+,\cdot),\\
x & \mapsto & (\pi_1(x),\pi_2(x),\ldots,\pi_t(x)).
\end{array}
\end{equation}
Let   $\pi_i^\ast$ (resp. $\pi^\ast$) denote the restriction of $\pi_i$ (resp. $\pi$) on $\mathbb{Z}_{h}^\ast$ for each $i\in[t]$.
Then there is an  isomorphism  as follows (cf. \cite{Ireland}):
\begin{equation}\label{equa3}
\begin{array}{rcl}
\pi^\ast:(\mathbb{Z}_{h}^\ast,\cdot) &\rightarrow & (\mathbb{Z}_{p_1^{s_1}}^\ast\times \mathbb{Z}_{p_2^{s_2}}^\ast
 \times\cdots\times\mathbb{Z}_{p_t^{s_t}}^\ast,\cdot),\\
x & \mapsto & (\pi_1^\ast(x),\pi_2^\ast(x),\ldots,\pi_t^\ast(x)).
\end{array}
\end{equation}

For a subset $S$ of $\mathbb{Z}_{h}$, let $S^{m\times n}$ be the set of all $m\times n$ matrices over $S$, and $S^{n}=S^{1\times n}$.
A matrix in $S^n$ is also called an $n$-dimensional row vector  over $S$.
 Let $I_r$ ($I$ for short) be the $r\times r$ identity matrix, $0_{m,n}$ ($0$ for short) the $m\times n$ zero matrix and $0_n=0_{n,n}$,
  and $\det(X)$ the determinant of a square matrix $X$ over $\mathbb{Z}_{h}$. Let $m\leq n$ and $a_i\in\mathbb{Z}_{h}$ for each $i\in[m]$. Define
$$\hbox{diag}(a_1,a_2,\ldots, a_m)_{m,n}
:=\left(\begin{array}{cccc|c}
a_1 &&&&0_{1,n-m}\\
& a_2&&&0_{1,n-m}\\
&&\ddots&&\vdots\\
&&&a_m&0_{1,n-m}
\end{array}\right)$$
and
$$\hbox{diag}(a_1,a_2,\ldots, a_m)_{n,m}
:=\left(\begin{array}{cccc}
a_1 &&&\\
& a_2&&\\
&&\ddots&\\
&&&a_m  \\ \hline
0_{n-m,1} &0_{n-m,1} &\cdots & 0_{n-m,1}
\end{array}\right).$$
 The set of $n\times n$ invertible matrices forms a group under matrix multiplication, called the {\it general linear group} of degree $n$ over $\mathbb{Z}_{h}$ and denoted by $G\!L_n(\mathbb{Z}_{h})$.

For $A\in\mathbb{Z}_{h}^{m\times n}$ and $B\in\mathbb{Z}_{h}^{n\times m}$, if $AB=I_m$, we say that $A$ has
a {\it right inverse} and $B$ is a right inverse of $A$. Similarly, if $AB=I_m$, then $B$ has a {\it left
inverse} and $A$ is a left inverse of $B$.
$\mathbb{Z}_{h}^{n}$ is called the $n$-dimensional row vector space over $\mathbb{Z}_{h}$. Let $\alpha_i\in\mathbb{Z}_{h}^{n}$ for each $i\in[m]$.
The vector subset $\{\alpha_1,\alpha_2,\ldots,\alpha_m\}$ is called {\it unimodular} if the matrix
$(\alpha_1^t, \alpha_2^t, \ldots, \alpha_m^t)^t$
has a right inverse, where $X^t$ is the transpose of $X$. By Lemma~\ref{lemma2.12} below, a matrix  $A\in\mathbb{Z}_{h}^{m\times n}$
has a right inverse if and only if all row vectors of $A$ are linearly independent in $\mathbb{Z}_{h}^{n}$.

Let $V\subseteq\mathbb{Z}_{h}^{n}$ be a {\it linear subset} ($\mathbb{Z}_{h}$-module). A {\it largest unimodular vector subset}
of $V$ is a unimodular vector subset of $V$ which has maximum number of vectors. The
dimension of $V$, denoted by $\dim(V)$, is the number of vectors in a largest unimodular
vector subset of $V$. Clearly, $\dim(V)=0$ if and only if $V$ does not contain a unimodular
vector. If a linear subset $X$ of $\mathbb{Z}_{h}^{n}$  has a unimodular basis with $m$ vectors, then $X$ is
called an $m$-{\it dimensional vector subspace} ($m$-{\it subspace} for short) of $\mathbb{Z}_{h}^{n}$. Every
$m$-subspace of $\mathbb{Z}_{h}^{n}$ is isomorphic to $\mathbb{Z}_{h}^{m}$. We define the $0$-subspace to be $\{0\}$.

The Erd\H{o}s-Ko-Rado theorem \cite{Erdos,Wilson}
 is a classical result in extremal set theory which obtained an upper bound on the cardinality of a family of $m$-subsets of a set such that every pairwise intersection has cardinality at least $r$ and describes exactly which families meet this bound. The results on Erd\H{o}s-Ko-Rado theorem have inspired much research \cite{Frankl,Godsi2,Huang-T,Tanaka,Vanhove}.

\begin{defin}
 Let $0\leq r\leq m\leq n$ and ${\mathbb{Z}_{h}^{n}\brack m}$ be the set of all $m$-subspaces of $\mathbb{Z}_{h}^{n}$.
 A family ${\cal F}\subseteq{\mathbb{Z}_{h}^{n}\brack m}$ is called $r$-{\it intersecting}
 if $\dim(A\cap B)\geq r$ for all $A,B\in{\cal F}$.
 \end{defin}

 When $t=1$ ($h$ is a prime power), Huang et al. \cite{Huang3} obtained an upper bound on the cardinality of an $r$-intersecting family in ${\mathbb{Z}_{h}^{n}\brack m}$ and described exactly which families meet this bound.

 \begin{thm}{\rm(See Theorem~1.1 in \cite{Huang3}.)}
 Let $\lfloor \frac{n}{2}\rfloor\geq m\geq r\geq 0$ and ${\cal F}\subseteq{\mathbb{Z}_{p^s}^{n}\brack m}$ be an $r$-intersecting family. Then
 $$|{\cal F}|\leq p^{(s-1)(n-m)(m-r)}{n-r\brack m-r}_{p},\;\hbox{where}\;
 {n-r\brack m-r}_{p}=\prod_{j=0}^{m-r-1}\frac{p^{n-r-j}-1}{p^{m-r-j}-1},$$
 and equality holds if and only if either {\rm(a)} ${\cal F}$ consists of all $m$-subspaces of $\mathbb{Z}_{p^s}^{n}$ which
contain a fixed $r$-subspace of $\mathbb{Z}_{p^s}^{n}$,  or {\rm(b)} $n=2m$ and ${\cal F}$ is the set of all $m$-subspaces of
$\mathbb{Z}_{p^s}^{n}$ contained in a fixed $(n-r)$-subspace of $\mathbb{Z}_{p^s}^{n}$.
 \end{thm}

For ${\cal I}\subseteq[t]$, by (\ref{equa2}), there exists the unique $x_{\cal I}\in\mathbb{Z}_h$ such that $\pi_i(x_{\cal I})=1$ if $i\in {\cal I}$,
 and $\pi_i(x_{\cal I})=0$ if $i\in[t]\setminus {\cal I}$.
 Let $0\leq2r\leq2m=n$ and ${\cal I}\subseteq[t]$.  Define
\begin{equation}\label{equa4}
{\cal F}_{(\alpha_1,\alpha_2,\ldots,\alpha_t)}^{(r,m,n,{\cal I})}
=\left\{\left(\begin{array}{cc}
X & 0_{m-r,r}\\
Y & x_{{\cal I}}I_{r}
\end{array}\right)\in{\mathbb{Z}_h^n\brack m} :
X\in\mathbb{Z}_h^{(m-r)\times(n-r)},Y\in J_{(\alpha_1,\alpha_2,\ldots,\alpha_t)}^{r\times(n-r)} \right\},
\end{equation}
where $\pi_i(x_{\cal I})=1$ and $\alpha_i=s_i$ if $i\in {\cal I}$;
 and $\pi_i(x_{\cal I})=0$ and $\alpha_i=0$ if $i\in[t]\setminus {\cal I}$.

As a generalization of Theorem~1.1, we give the Erd\H{o}so-Rado theorem for $\mathbb{Z}_{h}^{n}$ as follows.

 \begin{thm}\label{thm1.1}
 Let $\lfloor \frac{n}{2}\rfloor\geq m\geq r\geq 0$ and ${\cal F}\subseteq{\mathbb{Z}_{h}^{n}\brack m}$ be an $r$-intersecting family. Then
 $$|{\cal F}|\leq \prod_{i=1}^tp_i^{(s_i-1)(n-m)(m-r)}{n-r\brack m-r}_{p_i},$$
 and equality holds if and only if either {\rm(a)} ${\cal F}$ consists of all $m$-subspaces of $\mathbb{Z}_{h}^{n}$ which
contain a fixed $r$-subspace of $\mathbb{Z}_{h}^{n}$,  {\rm(b)} $n=2m$ and ${\cal F}$ is the set of all $m$-subspaces of
$\mathbb{Z}_{h}^{n}$ contained in a fixed $(n-r)$-subspace of $\mathbb{Z}_{h}^{n}$, or {\rm(c)} $n=2m$ and  there exists some $T\in G\!L_n(\mathbb{Z}_h)$ 
such that ${\cal F}=\left\{FT: F\in{\cal F}_{(\alpha_1,\alpha_2,\ldots,\alpha_t)}^{(r,m,n,{\cal I})}\right\}$ with ${\cal I}\not=\emptyset,[t]$, where ${\cal F}_{(\alpha_1,\alpha_2,\ldots,\alpha_t)}^{(r,m,n,{\cal I})}$
is given by (\ref{equa4}).
 \end{thm}

To prove Theorem~\ref{thm1.1}, we need to discuss generalized Grassmann graphs over $\mathbb{Z}_{h}$.
The Grassmann graph over a finite field plays an important role in geometry and combinatorics,  see \cite{Brouwer,Pankov,Wan2}. As a natural extension, the {\it generalized Grassmann graph} over  $\mathbb{Z}_{h}$, denoted by $G_r(m,n,\mathbb{Z}_{h})$, has the vertex set ${\mathbb{Z}_{h}^n\brack m}$, and two distinct vertices
are adjacent if their intersection is of dimension $>m-r$,
where  $2\leq r\leq m+1\leq n$.  Note that $G_2(m,n,\mathbb{Z}_{h})$ is the Grassmann graph $G(m,n,\mathbb{Z}_{h})$.

Let $V(\Gamma)$ denote the vertex set of a graph $\Gamma$. For $A,B\in V(\Gamma)$, we write $A\sim B$ if vertices $A$ and $B$ are adjacent. A {\it clique} of a graph $\Gamma$ is a complete subgraph of $\Gamma$. A clique ${\cal C}$ is {\it maximal} if there is no clique of $\Gamma$ which properly contains ${\cal C}$ as a subset. A {\it maximum clique} of $\Gamma$ is a clique of $\Gamma$ which has maximum cardinality. The {\it clique number} $\omega(\Gamma)$ of $\Gamma$ is the number of vertices in a maximum clique.

In the rest of this paper, we always assume that $m\leq n$ and $h$ is as in (\ref{equa1}).
The paper is organized as follows. In Section~2, we introduce the basic properties of
matrices over $\mathbb{Z}_{h}$ for later reference. In Section~3, we obtain the dimension
formula and the Anzahl theorem for subspaces of $\mathbb{Z}_{h}^n$. In Section~4, we calculate the clique number and determine geometric structures
of maximum cliques of $G_r(m,n,\mathbb{Z}_{h})$.   As a result,  Theorem~\ref{thm1.1} is proved.

\section{Matrices  over $\mathbb{Z}_h$}
For $A=(a_{uv})\in\mathbb{Z}_{h}^{m\times n}$, let
$\pi_i(A)=(\pi_i(a_{uv}))$ for each $i\in[t]$. Then $\pi_i$ is the natural surjective homomorphism from
$(\mathbb{Z}_{h}^{m\times n},+)$ to $(\mathbb{Z}_{p_i^{s_i}}^{m\times n},+)$ for each $i\in[t]$, and there is an isomorphism  as follows (cf. \cite{Brown}):
\begin{equation}\label{equa5}
\begin{array}{rcl}
\pi:(\mathbb{Z}_{h}^{m\times n},+) &\rightarrow & (\mathbb{Z}_{p_1^{s_1}}^{m\times n}\oplus \mathbb{Z}_{p_2^{s_2}}^{m\times n}
 \oplus\cdots\oplus\mathbb{Z}_{p_t^{s_t}}^{m\times n},+),\\
A & \mapsto & (\pi_1(A),\pi_2(A),\ldots,\pi_t(A)).
\end{array}
 \end{equation}
Let   $h_i=h/p_i^{s_i}$ and $\theta_i$  be the natural surjective homomorphism from
$(\mathbb{Z}_{h},+,\cdot)$ to $(\mathbb{Z}_{h_i},+,\cdot)$ for each $i\in[t]$. For $A=(a_{uv})\in\mathbb{Z}_{h}^{m\times n}$, let
$\theta_i(A)=(\theta_i(a_{uv}))$ for each $i\in[t]$. Then $\theta_i$ is the natural surjective homomorphism from
$(\mathbb{Z}_{h}^{m\times n},+)$ to $(\mathbb{Z}_{h_i}^{m\times n},+)$ for each $i\in[t]$, and
there is an isomorphism  as follows:
\begin{equation}\label{equa6}
\begin{array}{rcl}
(\pi_i,\theta_i):(\mathbb{Z}_{h}^{m\times n},+) &\rightarrow & (\mathbb{Z}_{p_i^{s_i}}^{m\times n}\oplus \mathbb{Z}_{h_i}^{m\times n},+),\\
A & \mapsto & (\pi_i(A),\theta_i(A)).
\end{array}
 \end{equation}
 For matrices $A\in\mathbb{Z}_{h}^{m\times n}$ and $B\in\mathbb{Z}_{h}^{n\times k}$, it is easy to see that
 \begin{equation}\label{equa7}
 \pi_i(AB)=\pi_i(A)\pi_i(B)\quad\hbox{and}\quad\theta_i(AB)=\theta_i(A)\theta_i(B)\;\hbox{for each}\;i\in[t].
 \end{equation}

\begin{lemma}\label{lemma2.1}{\rm(See Lemma~2.2 in \cite{Guo3}.)}
Every non-zero element $x$ in $\mathbb{Z}_{h}$ can be written as $x=u\prod_{i=1}^tp_i^{\alpha_i}$, where $u$ is a unit, and
$0\leq\alpha_i\leq s_i$ for each $i\in[t]$.
Moreover, the vector $(\alpha_1,\alpha_2,\ldots,\alpha_t)$ is unique
and $u$ is unique modulo the ideal $J_{(s_1-\alpha_1,s_2-\alpha_2,\ldots,s_t-\alpha_t)}$.
\end{lemma}

\begin{lemma}\label{lemma2.2}{\rm(See Theorem~2.7 in \cite{Guo3}.)}
Every matrix $A\in\mathbb{Z}_{h}^{m\times n}$ is equivalent to
\begin{equation}\label{equa8}
D={\rm diag}\,\left(\prod_{i=1}^tp_i^{\alpha_{i1}},\prod_{i=1}^tp_i^{\alpha_{i2}},\ldots,\prod_{i=1}^tp_i^{\alpha_{im}}\right)_{m,n},
\end{equation}
where $0\leq\alpha_{i1}\leq\alpha_{i2}\leq\cdots\leq\alpha_{im}\leq s_i$ for all $i\in[t]$.
Moreover, the array
$$
\Omega=((\alpha_{11},\alpha_{12},\ldots,\alpha_{1m}),(\alpha_{21},\alpha_{22},\ldots,\alpha_{2m}),
\ldots,(\alpha_{t1},\alpha_{t2},\ldots,\alpha_{tm}))
$$
is uniquely determined by $A$.
\end{lemma}

Let $A\in \mathbb{Z}_{h}^{m\times n}$ be a non-zero matrix. By Cohn's definition \cite{Cohn2}, the {\it inner rank} of $A$, denoted by $\rho(A)$,
 is the least integer $r$ such that
 \begin{equation}\label{equa10}
A = BC,\quad \hbox{where}\; B \in \mathbb{Z}_{h}^{m\times r} \;\hbox{and}\; C \in \mathbb{Z}_{h}^{r\times n}.
\end{equation}
Let $\rho(0) =0$. Any factorization as (\ref{equa10}) with $r=\rho(A)$ is called a {\it minimal factorization} of $A$.
For $A\in \mathbb{Z}_{h}^{m\times n}$, it is obvious that $\rho(A)\leq\min\{m,n\}$ and $\rho(A)=0$ if and only if $A=0$.

For matrices over $\mathbb{Z}_{h}$, from \cite{Cohn,Cohn2}, we deduce that
the following hold:
 \begin{equation}\label{equa11}
 \rho(A)=\rho(SAT),\;\hbox{where}\;S\;\hbox {and}\;T\;\hbox{are invertible matrices}.
 \end{equation}
 \begin{equation}\label{equa12}
   \rho(AB) \leq\min\{\rho(A),\rho(B)\}.
\end{equation}
$$
   \rho\left(\begin{array}{cc}
   A_{11} & A_{12}\\
  A_{21}  & A_{22}
  \end{array}\right) \geq\max\{\rho(A_{ij}) :1\leq i,j\leq 2\}.
$$

\begin{lemma}\label{lemma2.3}{\rm(See Lemma~2.8 in \cite{Guo3}.)}
Let $A=SDT\in\mathbb{Z}_{h}^{m\times n}$, where $S$ and $T$ are invertible, and $D$ is as in (\ref{equa8}).
Then the inner rank of $A$ is $\max\{c: (\alpha_{1c},\alpha_{2c},\ldots,\alpha_{tc})\not=(s_1,s_2,\ldots,s_t)\}$.
\end{lemma}

\begin{lemma}\label{lemma2.4}{\rm(See Lemma~2.9 in \cite{Guo3}.)}
Let $A\in\mathbb{Z}_{h}^{m\times n}$. Then
$$\rho(A)=\max\{\rho(\pi_i(A)): i=1,2,\ldots,t\}=\max\{\rho(\theta_i(A)): i=1,2,\ldots,t\}.$$
\end{lemma}

Let $A\in\mathbb{Z}_{h}^{m\times n}$. Denote by $I_k(A)$ the ideal in $\mathbb{Z}_{h}$
generated by all $k\times k$ minors of $A$.
Let $\hbox{Ann}_{\mathbb{Z}_{h}}(I_k(A)) =\{x\in \mathbb{Z}_{h}:xI_k(A)=0\}$ denote the annihilator of $I_k(A)$. The {\it McCoy rank} of $A$, denoted by $\hbox{rk}(A)$, is the following integer:
$$\hbox{rk}(A) = \max\{k : \hbox{Ann}_{\mathbb{Z}_{h}}(I_k(A)) =(0)\} .$$
Note that $\hbox{rk}(A) =\hbox{rk}(A^t)$; $\hbox{rk}(A) =\hbox{rk}(SAT)$ where $S$ and $T$ are invertible matrices; and $\hbox{rk}(A) =0$ if and only if $\hbox{Ann}_{\mathbb{Z}_{h}}(I_1(A)) \not=(0)$, see \cite{Brown}.

\begin{lemma}\label{lemma2.5}
Let $A=SDT\in\mathbb{Z}_{h}^{m\times n}$, where $S$ and $T$ are invertible, and $D$ is as in (\ref{equa8}).
Then the McCoy rank of $A$ is  $\max\{c: (\alpha_{1c},\alpha_{2c},\ldots,\alpha_{tc})=(0,0,\ldots,0)\}.$
\end{lemma}
\proof  Let $\ell=\max\{c:(\alpha_{1c},\alpha_{2c},\ldots,\alpha_{tc})=(0,0,\ldots,0)\}$.
Then $I_\ell(A)=\mathbb{Z}_{h}$ and $I_{\ell+1}(A)=(\prod_{i=1}^tp_i^{\alpha_{i,\ell+1}})$, which imply that $\hbox{Ann}_{\mathbb{Z}_{h}}(I_\ell(A))=(0)$ and $\hbox{Ann}_{\mathbb{Z}_{h}}(I_{\ell+1}(A))\not=(0)$. Similarly, $\hbox{Ann}_{\mathbb{Z}_{h}}(I_{k}(A))\not=(0)$ for $\ell+1 \leq k \leq m$. It follows that ${\rm rk}(A)=\ell$.$\qed$

For matrices $A, B$ over $\mathbb{Z}_{h}$, by Lemmas~\ref{lemma2.2},~\ref{lemma2.3}  and~\ref{lemma2.5}, if $\rho(A)={\rm rk}(A)$, then
$$
   \rho\left(\begin{array}{cc}
   A & 0\\ 0& B
   \end{array}\right)=\rho(A)+\rho(B)\quad\hbox{and}\quad
   \hbox{rk}\left(\begin{array}{cc}
   A &0 \\ 0 & B
   \end{array}\right)=\hbox{rk}(A)+\hbox{rk}(B).
$$

\begin{lemma}\label{lemma2.6}{\rm(See \cite{Rosen}.)}
Let  $\pi_i^\ast$ (resp. $\pi^\ast$) also denote the restriction of $\pi_i$ (resp. $\pi$) on $G\!L_{n}(\mathbb{Z}_{h})$  for each $i\in[t]$.
Then $\pi_i^\ast$ is a natural surjective homomorphism  from
$(G\!L_{n}(\mathbb{Z}_{h}),\cdot)$ to $(G\!L_{n}(\mathbb{Z}_{p_i^{s_i}}),\cdot)$ for each $i\in[t]$, and
$\pi^\ast$ is an isomorphism from $(G\!L_{n}(\mathbb{Z}_{h}),\cdot)$ to $(G\!L_{n}(\mathbb{Z}_{p_1^{s_1}})\times G\!L_{n}(\mathbb{Z}_{p_2^{s_2}})\times\cdots\times G\!L_{n}(\mathbb{Z}_{p_t^{s_t}}),\cdot)$.
\end{lemma}

\begin{lemma}\label{lemma2.7}{\rm(See Corollary~2.21 in \cite{Brown}.)}
Let $A\in \mathbb{Z}_h^{n\times n}$. Then $A\in G\!L_n(\mathbb{Z}_{h})$ if and only if $\det(A)\in\mathbb{Z}_h^\ast$.
\end{lemma}

\begin{lemma}\label{lemma2.8}
Let $A\in\mathbb{Z}_{h}^{m\times n}$. Then
$${\rm rk}(A)=\min\{{\rm rk}(\pi_i(A)): i=1,2,\ldots,t\}=\min\{{\rm rk}(\theta_i(A)): i=1,2,\ldots,t\}.$$
\end{lemma}
\proof
By Lemma~\ref{lemma2.2}, we may assume that $A=SDT$, where $S$ and $T$ are invertible, and $D$ is as in (\ref{equa8}). Let ${\rm rk}(A)=\ell$.
For each $i\in[t]$, we write $w_{ic}=\pi_i(\prod_{j=1}^tp_j^{\alpha_{jc}})$ for all $c=\ell+1,\ldots,m$.
By Lemma~\ref{lemma2.6} and~(\ref{equa7}), $\pi_{i}(S)=\pi_{i}^\ast(S)\in G\!L_m(\mathbb{Z}_{p_i^{s_i}}),\pi_{i}(T)=\pi_{i}^\ast(T)\in G\!L_n(\mathbb{Z}_{p_i^{s_i}})$, and
$$
\pi_i(A)=\pi_{i}(S)\hbox{diag}(1,\ldots,1,w_{i,\ell+1},\dots,w_{im})_{m,n}\pi_{i}(T).
$$
By Lemma~\ref{lemma2.5}, we obtain $${\rm rk}(\pi_i(A))={\rm rk}(\hbox{diag}(1,\ldots,1,w_{i,\ell+1},\dots,w_{im})_{m,n})\geq \ell={\rm rk}(A).$$
Therefore, ${\rm rk}(A)\leq\min\{{\rm rk}(\pi_i(A)): i=1,2,\ldots,t\}$.
Since $(\alpha_{1,\ell+1},\alpha_{2,\ell+1},\ldots,\alpha_{t,\ell+1})\not=(0,0,\ldots,0)$,
there exists some $i_0\in[t]$ such that $\alpha_{i_0,\ell+1}\not=0$, which implies that $w_{i_0,\ell+1}\not\in\mathbb{Z}_{p_{i_0}^{s_{i_0}}}^\ast$.
It follows that ${\rm rk}(\pi_{i_0}(A))=\ell={\rm rk}(A)$.

For each $i\in[t]$, we write $x_{ic}=\theta_i(\prod_{j=1}^tp_j^{\alpha_{jc}})$ for all $c=\ell+1,\ldots,m$.
By Lemma~\ref{lemma2.7} and~(\ref{equa7}), $\theta_{i}(S)\in G\!L_m(\mathbb{Z}_{h_i}),\theta_{i}(T)\in G\!L_n(\mathbb{Z}_{h_i})$, and
$$
\theta_i(A)=\theta_{i}(S)\hbox{diag}(1,\ldots,1,x_{i,\ell+1},\dots,x_{im})_{m,n}\theta_{i}(T),
$$
where  $h_i=h/p_i^{s_i}$.
By Lemma~\ref{lemma2.5} again, we obtain $${\rm rk}(\theta_i(A))={\rm rk}(\hbox{diag}(1,\ldots,1,x_{i,\ell+1},\dots,x_{im})_{m,n})\geq \ell={\rm rk}(A).$$
Therefore, ${\rm rk}(A)\leq\min\{{\rm rk}(\theta_i(A)): i=1,2,\ldots,t\}$.
Since $(\alpha_{1,\ell+1},\alpha_{2,\ell+1},\ldots,\alpha_{t,\ell+1})\not=(0,0,\ldots,0)$,
there exists some $i_0\in[i]$ such that $\alpha_{i_0,\ell+1}\not=0$, which implies that $x_{i,\ell+1}\not\in\mathbb{Z}_{h_i}^\ast$ for all $i\not=i_0$. It follows that
${\rm rk}(\theta_i(A))=\ell={\rm rk}(A)$ for all $i\not=i_0$.
Hence, the desired result follows. $\qed$

\begin{lemma}\label{lemma2.9}{\rm(See Lemma~2.4 in \cite{Huang3}.)}
Let $X\in\mathbb{Z}_{p^s}^{m\times n}$. Then $X$ has a right inverse if and only if ${\rm rk}(X)=m$.
\end{lemma}

\begin{lemma}\label{lemma2.10}
Let $X\in\mathbb{Z}_{h}^{m\times n}$. Then the following hold:
\begin{itemize}
\item[\rm(i)]
$X$ has a right inverse if and only if ${\rm rk}(X)=m$.

 \item[\rm(ii)]
 $X$ has a right inverse if and only if $\pi_i(X)$ has a right inverse for all $i\in[t]$.

\item[\rm(iii)]
When $m=n$, $X$ is invertible if and only if $\pi_i(X)$ is invertible for all $i\in[t]$. Moreover, if $X$ is
invertible, then $\pi_i(X^{-1})=(\pi_i(X))^{-1}$ for each $i\in[t]$.
\end{itemize}
\end{lemma}
\proof (i).
If $X$ has a right inverse, then (\ref{equa7}) implies that $\pi_i(X)$ has a right inverse for all $i\in[t]$.
By Lemma~\ref{lemma2.9},  ${\rm rk}(\pi_i(X))=m$ for all $i\in[t]$. By Lemma~\ref{lemma2.8}, ${\rm rk}(X)=m$.

Conversely, suppose ${\rm rk}(X)=m$. By Lemma~\ref{lemma2.2}, there are $S\in G\!L_m(\mathbb{Z}_h)$ and $T\in G\!L_n(\mathbb{Z}_h)$ such that $X=S(I_m,0_{m,n-m})T$.
Pick $Y=T^{-1}\left(\begin{array}{c}
I_{m}\\
0_{n-m,m}
\end{array}\right)S^{-1}$.
Then $XY=I_m$, which implies that $X$ has a right inverse.

(ii). If $X$ has a right inverse, then  $\pi_i(X)$ has a right inverse for all $i\in[t]$.
Conversely, suppose that $\pi_i(X)$ has a right inverse for all $i\in[t]$. By the proof of (i), ${\rm rk}(X)=m$, which implies that
$X$ has a right inverse by (i).

(iii). By Lemma~\ref{lemma2.6}, $X$ is invertible if and only if $\pi_i(X)$ is invertible for all $i\in[t]$.
Suppose that $X\in\mathbb{Z}_{h}^{n\times n}$ is invertible. Then $XX^{-1}=I_n$. From (\ref{equa7}), we deduce  that
$\pi_i(XX^{-1})=\pi_i(X)\pi_i(X^{-1})=I_n$ for each $i\in[t]$. Therefore, we have  $\pi_i(X^{-1})=(\pi_i(X))^{-1}$ for each $i\in[t]$. $\qed$

\begin{lemma}\label{lemma2.11}
Let $X\in\mathbb{Z}_{h}^{m\times n}$. Then the following hold:
\begin{itemize}
\item[\rm(i)]
 $X$ has a right inverse if and only if $\theta_i(X)$ has a right inverse for all $i\in[t]$.

\item[\rm(ii)]
When $m=n$, $X$ is invertible if and only if $\theta_i(X)$ is invertible for all $i\in[t]$. Moreover, if $X$ is
invertible, then $\theta_i(X^{-1})=(\theta_i(X))^{-1}$ for each $i\in[t]$.
\end{itemize}
\end{lemma}
\proof (i). If $X$ has a right inverse, then (\ref{equa7}) implies that $\theta_i(X)$ has a right inverse for all $i\in[t]$.
Conversely, suppose that $\theta_i(X)$ has a right inverse for all $i\in[t]$. By Lemma~\ref{lemma2.10} (i), ${\rm rk}(\theta_i(X))=m$ for all $i\in[t]$,
which imply that ${\rm rk}(X)=m$ by Lemma~\ref{lemma2.8}. By Lemma~\ref{lemma2.10} (i) again, $X$ has a right inverse.

(ii). Similarly, by  Lemma~\ref{lemma2.7}, we may prove that $X$ is invertible if and only if $\theta_i(X)$ is invertible for all $i\in[t]$.
Similar to the proof of Lemma~\ref{lemma2.10} (iii), we have $\theta_i(X^{-1})=(\theta_i(X))^{-1}$ for each $i\in[t]$. $\qed$

\begin{lemma}\label{lemma2.12}
Let $\alpha_1,\alpha_2,\ldots,\alpha_m\in\mathbb{Z}_{h}^{n}$. Then the following hold:
\begin{itemize}
\item[\rm(i)]
The matrix  $A=(\alpha_1^t, \alpha_2^t, \ldots, \alpha_m^t)^t$
has a right inverse if and only if $\alpha_1,\alpha_2,\ldots,\alpha_m$ are linearly independent.

\item[\rm(ii)]
If $\alpha_1,\alpha_2,\ldots,\alpha_m$ are linearly independent, then $\alpha_1,\alpha_2,\ldots,\alpha_m$ can be extended to a basis of $\mathbb{Z}_{h}^{n}$.
\end{itemize}
\end{lemma}
\proof (i). If $A$ has a right inverse, then there exists some $B\in\mathbb{Z}_h^{n\times m}$ such that $AB=I_m$. From $\sum_{i=1}^mx_i\alpha_i=0$, we deduce that $(x_1,x_2,\ldots,x_m)A=0$, which implies that $(x_1,x_2,\ldots,x_m)=(x_1,x_2,\ldots,x_m)AB=0$, and therefore $\alpha_1,\alpha_2,\ldots,\alpha_m$ are linearly independent.

Conversely, suppose that $\alpha_1,\alpha_2,\ldots,\alpha_m$ are linearly independent. Without loss of generality, by Lemma~\ref{lemma2.2}, we may assume that $A=SDT$, where $S$ and $T$ are invertible, and $D$ is as in (\ref{equa8}). By Lemma~\ref{lemma2.5}, ${\rm rk}(A)=m$ if and only if $(\alpha_{1m},\alpha_{2m},\ldots,\alpha_{tm})=(0,0,\ldots,0)$,
 which implies that $A$ has a right inverse if and only if $(\alpha_{1m},\alpha_{2m},\ldots,\alpha_{tm})=(0,0,\ldots,0)$ by Lemma~\ref{lemma2.9}. Assume that $(\alpha_{1m},\alpha_{2m},\ldots,\alpha_{tm})\not=(0,0,\ldots,0)$, then $x_m=\prod_{i=1}^tp_i^{s_i-\alpha_{im}}\not=0$ such that $(0,\ldots,0,x_m)S^{-1}\not=0$ and $(0,\ldots,0,x_m)S^{-1}A=0$,
which implies that $\alpha_1,\alpha_2,\ldots,\alpha_m$ are not linearly independent, a contradiction.

(ii). Since $\alpha_1,\alpha_2,\ldots,\alpha_m$ are linearly independent, $A=(\alpha_1^t, \alpha_2^t, \ldots, \alpha_m^t)^t$ has a right inverse,
 which implies that ${\rm rk}(A)=m$  by Lemma~\ref{lemma2.10} (i).
By Lemma~\ref{lemma2.2}, there are $S\in G\!L_m(\mathbb{Z}_h)$ and $T\in G\!L_n(\mathbb{Z}_h)$ such that $A=S(I_m,0_{m,n-m})T$.
By Lemma~\ref{lemma2.7}, $$
\widetilde{A}=\left(\begin{array}{c}
A\\ (0_{n-m,m},I_{n-m})T
\end{array}\right)
=\left(\begin{array}{cc}
S & 0_{m,n-m}\\ 0_{n-m,m} & I_{n-m}
\end{array}\right)T$$
is invertible, which implies that all row vectors of $\widetilde{A}$
are a basis of $\mathbb{Z}_{h}^{n}$, and therefore $\alpha_1,\alpha_2,\ldots,\alpha_m$ can be extended to a basis of $\mathbb{Z}_{h}^{n}$. $\qed$

\section{Subspaces of $\mathbb{Z}_{h}^{n}$}
In this section, we study subspaces of $\mathbb{Z}_{h}^{n}$ and obtain some useful results.
Notations and terminologies will be adopted from \cite{Huang3}.

Let $\alpha_1,\alpha_2,\ldots,\alpha_m$ be $m$ linearly independent row vectors in $\mathbb{Z}_{h}^{n}$,
and $[\alpha_1,\alpha_2,\ldots,\alpha_m]$ be the $\mathbb{Z}_{h}$-module generated by $\alpha_1,\alpha_2,\ldots,\alpha_m$.
Then $X=[\alpha_1,\alpha_2,\ldots,\alpha_m]$  is an $m$-subspace of $\mathbb{Z}_{h}^{n}$, and the matrix
$(\alpha_1^t, \alpha_2^t, \ldots, \alpha_m^t)^t$
is called a {\it matrix representation} of $X$. We use an $m\times n$ matrix $A$ with
${\rm rk}(A)=m$ to represent an $m$-subspace of $\mathbb{Z}_{h}^{n}$. For $A,B\in\mathbb{Z}_{h}^{m\times n}$ with ${\rm rk}(A)={\rm rk}(B)=m$,
both $A$ and $B$ represent the same $m$-subspace if and only if there exists an $S\in G\!L_m(\mathbb{Z}_h)$ such that $B=SA$.
For convenience, if $A$ is a matrix representation of an $m$-subspace $X$, then we write $X=A$.

For an $m$-subspace $X$ of $\mathbb{Z}_{h}^{n}$,  by Lemma~\ref{lemma2.2}, there exists a $T\in G\!L_n(\mathbb{Z}_h)$
such that $X$ has matrix representations
$$X=(0_{m,n-m},I_m)T=S(0_{m,n-m},I_m)T\; \hbox{for all}\; S\in G\!L_m(\mathbb{Z}_h).$$
By Lemma~\ref{lemma2.1}, every $m$-subspace $X$ of $\mathbb{Z}_{h}^{n}$ has the
unique matrix representation $(A_1,I_m)T$, which is a {\it row-reduced echelon form}, where $T$ is a permutation matrix and $A_1\in\mathbb{Z}_{h}^{m\times(n-m)}$.

Let $X$ and $Y$ be two subspaces of $\mathbb{Z}_{h}^{n}$. A {\it join} of $X$ and $Y$ is a minimum dimensional subspace
containing $X$ and $Y$ in $\mathbb{Z}_{h}^{n}$. In general, the join of two subspaces in $\mathbb{Z}_{h}^{n}$ is not unique,
see \cite{Huang3} for examples.
Denoted by $X\vee Y$ the set of all joins of subspaces $X$ and $Y$ with the same minimum
dimension $\dim(X\vee Y )$. Then $X\vee Y=Y\vee X$, and $X\vee Y=\{Y\}$ if $X\subseteq Y$.
Clearly, $X\cap Y$ contains a subspace of dimension $\dim(X\cap Y)$. In general, $X\cap Y$ is a linear subset but it may not be a subspace,
see \cite{Huang3} for examples.

\begin{lemma}\label{lemma3.1} {\rm(See Lemma~3.2 in \cite{Huang3}.)}
Let $1\leq m\leq k<n$. Suppose that $A$ is a $k$-subspace of $\mathbb{Z}_{p^s}^{n}$ and $B$ is an
$m$-subspace of $\mathbb{Z}_{p^s}^{n}$ with $B\not\subseteq A$. Then there is a $T\in G\!L_n(\mathbb{Z}_{p^s})$ such that
$$
A=(0_{k,n-k},I_k)T\quad\hbox{and}\quad B=(D,I_m)T,
$$
where $D={\rm diag}\left(p^{\alpha_{1}},p^{\alpha_{2}},\ldots,p^{\alpha_{r}},0,\ldots,0\right)_{m,n-m},1\leq r\leq \min\{m,n-k\}$ and $0\leq\alpha_{1}\leq\alpha_{2}\leq\cdots\leq\alpha_{r}\leq \max\{s-1,1\}$.
\end{lemma}

\begin{lemma}\label{lemma3.2}
Let $1\leq m\leq k<n$. Suppose that $A$ is a $k$-subspace of $\mathbb{Z}_{h}^{n}$ and $B$ is an
$m$-subspace of $\mathbb{Z}_{h}^{n}$ with $B\not\subseteq A$. Then there is a $T\in G\!L_n(\mathbb{Z}_h)$ such that
\begin{equation}\label{equa15}
A=(0_{k,n-k},I_k)T\quad\hbox{and}\quad B=(D,I_m)T,
\end{equation}
where \begin{equation}\label{equa16}
D={\rm diag}\left(\prod_{i=1}^tp_i^{\alpha_{i1}},\prod_{i=1}^tp_i^{\alpha_{i2}},\ldots,\prod_{i=1}^tp_i^{\alpha_{ir}},0,\ldots,0\right)_{m,n-m}
\end{equation}
with $1\leq r\leq \min\{m,n-k\},\rho(D)=r,$ and
 $0\leq\alpha_{i1}\leq\alpha_{i2}\leq\cdots\leq\alpha_{ir}\leq s_i$ for all $i\in[t]$.
\end{lemma}
\proof For each $i\in[t]$, by Lemma~\ref{lemma2.10}, $\pi_i(A)$ is a $k$-subspace of $\mathbb{Z}_{p_i^{s_i}}^{n}$ and $\pi_i(B)$
is an $m$-subspace of $\mathbb{Z}_{p_i^{s_i}}^{n}$. By Lemma~\ref{lemma3.1},
there is a $T_i\in G\!L_n(\mathbb{Z}_{p_i^{s_i}})$ such that
$$
\pi_i(A)=(0_{k,n-k},I_k)T_i\quad\hbox{and}\quad \pi_i(B)=(D_i,I_m)T_i,
$$
where $D_i={\rm diag}\,\left(p_i^{\alpha_{i1}},p_i^{\alpha_{i2}},\ldots,p_i^{\alpha_{ir_i}},0,\ldots,0\right)_{m,n-m},1\leq r_i\leq \min\{m,n-k\}$ and $0\leq\alpha_{i1}\leq\alpha_{i2}\leq\cdots\leq\alpha_{ir_i}\leq \max\{s_i-1,1\}$.

Let $r=\max\{r_i:i=1,2,\ldots,t\}$ and $\alpha_{ij}=s_i$ for all $i\in[t]$ and $j=r_i+1,\ldots,r$. Then $1\leq r\leq \min\{m,n-k\}$ and
 $0\leq\alpha_{i1}\leq\alpha_{i2}\leq\cdots\leq\alpha_{ir}\leq s_i$ for all $i\in[t]$.
 For each  $c\in[r]$, by (\ref{equa2}), there exists the unique $d_c\in\mathbb{Z}_{h}$ such that $\pi_i(d_c)=p_i^{\alpha_{ic}}$ for all $i\in[t]$.
Let $\widetilde{D}={\rm diag}(d_1,d_2,\ldots,d_r,0,\ldots,0)_{m,n-m}$. Then $\pi_i(\widetilde{D})=D_i$ for all $i\in[t]$,
and $\rho(\widetilde{D})=r$ by Lemma~\ref{lemma2.4}.

For each  $c\in[r]$, let $\pi_i(\prod_{j=1}^tp_j^{\alpha_{jc}})=d_{ic}p_i^{\alpha_{ic}}$,
where $d_{ic}=\pi_i(\prod_{j\in[t]\setminus\{i\}}p_j^{\alpha_{jc}})\in\mathbb{Z}_{p_i^{s_i}}^\ast$
 for all $i\in[t]$. By (\ref{equa3}),
there exist $u_{c}\in\mathbb{Z}_h^\ast$ such that $\pi_i(u_c)=\pi_i^\ast(u_c)=d_{ic}$ for  all $i\in[t]$.
Since $\pi_i(u_cd_c)=\pi_i(\prod_{j=1}^tp_j^{\alpha_{jc}})$  for all $i\in[t]$, by (\ref{equa2}) again, we have $u_cd_c=\prod_{j=1}^tp_j^{\alpha_{jc}}$.

By Lemma~\ref{lemma2.6}, there exists the unique $\widetilde{T}\in G\!L_n(\mathbb{Z}_h)$ such that $\pi_i(\widetilde{T})=T_i$ for all $i\in[t]$.
By (\ref{equa7}), we have
$$\pi_i(A)=\pi_i((0_{k,n-k},I_k)\widetilde{T})\quad\hbox{and}\quad \pi_i(B)=\pi_i((\widetilde{D},I_m)\widetilde{T})\;\hbox{for all}\;i\in[t],$$
which imply that $\pi(A)=\pi((0_{k,n-k},I_k)\widetilde{T})$ and $\pi(B)=\pi((\widetilde{D},I_m)\widetilde{T}).$
From (\ref{equa5}),  we deduce that $A=(0_{k,n-k},I_k)\widetilde{T}$ and $B=(\widetilde{D},I_m)\widetilde{T}.$
Let $D={\rm diag}(u_1d_1,u_2d_2,\ldots,u_rd_r,0,\ldots,0)_{m,n-m}$ and
$T={\rm diag}\,(u_1^{-1},u_2^{-1},\ldots,u_r^{-1},1,\ldots,1)_{n,n}\widetilde{T}$. Then $D$ is as in (\ref{equa16}) such that
 $$A=(0_{k,n-k},I_k)T\quad\hbox{and}\quad B=(D,I_m)T.$$
Therefore, we complete the proof of this lemma. $\qed$

\begin{thm}\label{lemma3.3}{\rm(Dimensional formula.)}
Let $A$ and $B$ be two subspaces of $\mathbb{Z}_{h}^n$. Then
\begin{equation}\label{equa17}
\dim(A\vee B)=\dim(A) + \dim(B) - \dim(A\cap B) = \rho\left(\begin{array}{c}A\\ B \end{array}\right).
\end{equation}
Moreover,
\begin{eqnarray} \label{equa17-2}
\dim(A\cap B)&=&
\min\{\dim(\pi_i(A)\cap\pi_i(B)):i=1,2,\ldots,t\}, \nonumber \\
&=&\min\{\dim(\theta_i(A)\cap\theta_i(B)):i=1,2,\ldots,t\}.
\end{eqnarray}
\end{thm}
\proof Let $k=\dim(A)$ and $m=\dim(B)$. Without loss of generality, we may assume that $1\leq m\leq k<n$ and $B\not\subseteq  A$. By Lemma~\ref{lemma3.2}, there is a $T\in G\!L_n(\mathbb{Z}_h)$ such that $A$ and $B$ are as in (\ref{equa15}). Let
$D'={\rm diag}\left(\prod_{i=1}^tp_i^{\alpha_{i1}},\prod_{i=1}^tp_i^{\alpha_{i2}},\ldots,\prod_{i=1}^tp_i^{\alpha_{ir}}\right)_{r,r}$,
where $1\leq r\leq \min\{m,n-k\},\rho(D')=r,$ and
 $0\leq\alpha_{i1}\leq\alpha_{i2}\leq\cdots\leq\alpha_{ir}\leq s_i$ for all $i\in[t]$. Without loss of
generality, we may assume that $T=I_n$ and $m\geq2$. Hence we can assume further that
$A$ and $B$ have matrix representations
\begin{equation}\label{equa18}
A=(0_{k,n-k},I_k)\quad\hbox{and}\quad
B=\left(\begin{array}{cccc}
D' & 0_{r,n-m-r} & I_r & 0_{r,m-r}\\
0_{m-r,r} & 0_{m-r,n-m-r} & 0_{m-r,r} & I_{m-r}
\end{array}\right),
\end{equation}
respectively. Clearly, the $(m-r)$-subspace $(0_{m-r,n-m+r},I_{m-r})$ is contained in $A\cap B$, and
the $(k+r)$-subspace $\left(\begin{array}{ccc}
I_r & 0_{r,n-k-r}  & 0_{r,k}\\
0_{k,r} &  0_{k,n-k-r} & I_{k}
\end{array}\right)$
contains $A$ and $B$. It follows that $\dim(A\cap B)\geq m-r$ and $\dim(A\vee B)\leq k+r$.

Let $\alpha\in A\cap B$ be an $n$-dimensional row vector. Then there are matrices $C_1\in\mathbb{Z}_{h}^k$ and
$C_2=(C_{21},C_{22})\in\mathbb{Z}_{h}^m$, where $C_{21}\in\mathbb{Z}_{h}^r$ and $C_{22}\in\mathbb{Z}_{h}^{m-r}$,
 such that $\alpha=C_1A=(0_{1,n-k},C_1)$ and $\alpha=C_2B=(C_{21}D',0_{1,n-m-r},C_{21},C_{22})$. By $r\leq n-k$, we have  $C_{21}D'=0_{1,r}$ and
therefore $\alpha=(0_{1,n-m},C_{21},C_{22})$. Write $C_{21}=(y_1,\ldots,y_r)$.
From  $C_{21}D'=0_{1,r}$ and $\alpha_{i1}\leq\alpha_{i2}\leq\cdots\leq\alpha_{ir}$ for all $i\in[t]$, we deduce that $y_j\prod_{i=1}^tp_i^{\alpha_{ir}}=0$ for all $j\in[r]$. It follows that $y_j\in J_{(s_1-\alpha_{1r},s_2-\alpha_{2r},\ldots,s_t-\alpha_{tr})}$  for all $j\in[r]$.
For any $n$-dimensional row vector $\alpha\in A\cap B$, there are $y_1,\ldots,y_r\in J_{(s_1-\alpha_{1r},s_2-\alpha_{2r},\ldots,s_t-\alpha_{tr})}$ and $z_1,\ldots,z_{m-r}\in\mathbb{Z}_h$ such that
 $$\alpha=(0,\ldots,0,y_1,\ldots,y_r,z_1,\ldots,z_{m-r}).$$
Since $\rho(D')=r$, we have $(\alpha_{1r},\alpha_{2r},\ldots,\alpha_{tr})\not=(s_1,s_2,\ldots,s_t)$, which implies that $(s_1-\alpha_{1r},s_2-\alpha_{2r},\ldots,s_t-\alpha_{tr})\not=(0,0,\ldots,0)$. It follows that $\dim(A\cap B)\leq m-r$. Hence, we have $\dim(A\cap B)=m-r$.

Let $d=\dim(A\vee B)$. Then there exists a $d$-subspace $W\in A\vee B$. Since $B\not\subseteq A$,
one obtains $d>k$. So, $W$ has a matrix representation
$$W=\left(\begin{array}{c}
W_1\\ A
\end{array}\right)
=\left(\begin{array}{cc}
W_{11} & 0_{d-k,k}\\
0_{k,n-k} & I_k
\end{array}\right),$$
where $W_{11}\in\mathbb{Z}_{h}^{(d-k)\times(n-k)}$ has a right inverse. By $B\subseteq W$ and (\ref{equa18}), the $r$-subspace
$(D',0_{r,n-m-r},I_r,0_{r,m-r})\subseteq W$. Therefore, there is a matrix $Q=(Q_1,Q_2)\in\mathbb{Z}_{h}^{r\times d}$, where
$Q_1\in\mathbb{Z}_{h}^{r\times(d-k)}$ and $Q_2\in\mathbb{Z}_{h}^{r\times k}$, such that $(D',0_{r,n-m-r},I_r,0_{r,m-r})=QW=(Q_1W_{11},Q_2)$.
Hence $(D',0_{r,n-k-r})=Q_1W_{11}$. By (\ref{equa12}), $r=\rho(D')=\rho(Q_1W_{11})\leq\rho(W_{11})\leq d-k$, which implies that
$d=\dim(A\vee B)\geq k+r$. So, $\dim(A\vee B)=k+r$.

From (\ref{equa11}) and (\ref{equa18}), we deduce that $$\rho\left(\begin{array}{c}A\\ B \end{array}\right)
=\rho\left(\begin{array}{ccc}
I_k&0&0\\
0&D_1 &  0\\
0 & 0 & 0
\end{array}\right)=k+r.$$
Therefore, we complete the proof of  (\ref{equa17}).

Let us continue to assume that $A$ and $B$ have matrix representations as in (\ref{equa18}).
From $r=\rho(D')$, we deduce that $\dim(A\cap B)=m-\rho(D')$.
By Lemma~\ref{lemma2.10}, $\pi_i(A)$ and $\pi_i(B)$ are two subspaces of $\mathbb{Z}_{p_i^{s_i}}^n$ for all $i\in[t]$.
Similarly, we can prove that $\dim(\pi_i(A)\cap \pi_i(B))=m-\rho(\pi_i(D'))$ for all $i\in[t]$.
By Lemma~\ref{lemma2.8}, we have
\begin{eqnarray*}
\dim(A\cap B)&=&m-\rho(D')=m-\max\{\rho(\pi_i(D')): i=1,2,\ldots,t\}\\
&=&\min\{\dim(\pi_i(A)\cap\pi_i(B)):i=1,2,\ldots,t\}.
\end{eqnarray*}
Similarly, we have
$\dim(A\cap B)=\min\{\dim(\theta_i(A)\cap\theta_i(B)):i=1,2,\ldots,t\}$.
$\qed$

\begin{lemma}\label{lemma3.4}{\rm(See Theorem~3.5 in \cite{Huang3}.)}
Let $1\leq m<n$. Then the number of $m$-subspaces of $\mathbb{Z}_{p^s}^n$ is $p^{(s-1)m(n-m)}{n\brack m}_p$.
\end{lemma}

\begin{thm}\label{lemma3.5}
Let $1\leq m\leq k<n$. Then the following hold:
\begin{itemize}
\item[\rm(i)]
The number of $m$-subspaces of $\mathbb{Z}_{h}^n$ is $\prod_{i=1}^tp_i^{(s_i-1)m(n-m)}{n\brack m}_{p_i}.$

\item[\rm(ii)]
In $\mathbb{Z}_{h}^n$, the number of $m$-subspaces in a given $k$-subspace is $\prod_{i=1}^tp_i^{(s_i-1)m(k-m)}{k\brack m}_{p_i}.$

\item[\rm(iii)]
In $\mathbb{Z}_{h}^n$, the number of $k$-subspaces containing a fixed $m$-subspace is $$\prod_{i=1}^tp_i^{(s_i-1)(k-m)(n-k)}{n-m\brack k-m}_{p_i}.$$
\end{itemize}
\end{thm}
\proof
(i).  Suppose that $\pi$ is as in (\ref{equa5}). By Lemma~\ref{lemma2.10} and matrix representations of subspaces, $\pi$ induces a bijective map
as follows:
\begin{equation}\label{equa19}
\begin{array}{rcl}
\pi:{\mathbb{Z}_{h}^{n}\brack m} &\rightarrow & {\mathbb{Z}_{p_1^{s_1}}^{n}\brack m}\times {\mathbb{Z}_{p_2^{s_2}}^{n}\brack m}
 \times\cdots\times{\mathbb{Z}_{p_t^{s_t}}^{n}\brack m},\\
X & \mapsto & (\pi_1(X),\pi_2(X),\ldots,\pi_t(X)).
\end{array}
\end{equation}
By Lemma~\ref{lemma3.4}, the number of $m$-subspaces of $\mathbb{Z}_{p_i^{s_i}}^n$ is $p_i^{(s_i-1)m(n-m)}{n\brack m}_{p_i}$ for all $i\in[t]$, which imply that the number of $m$-subspaces of $\mathbb{Z}_{h}^n$ is
$\prod_{i=1}^tp_i^{(s_i-1)m(n-m)}{n\brack m}_{p_i}.$

(ii). Since every $k$-subspace of $\mathbb{Z}_{h}^n$ is isomorphic to $\mathbb{Z}_{h}^k$, the desired result follows by (i).

(iii). Let $${\cal S}=\left\{(A,B):A\in{\mathbb{Z}_{h}^{n}\brack m},B\in{\mathbb{Z}_{h}^{n}\brack k},A\subseteq B\right\}.$$
We compute the cardinality  of ${\cal S}$ in two ways. By (i) and (ii), the number of $k$-subspaces containing a fixed $m$-subspace is
$$\frac{|{\cal S}|}{\left|{\mathbb{Z}_{h}^{n}\brack m}\right|}=\frac{\left|{\mathbb{Z}_{h}^{n}\brack k}\right|\cdot\left|{\mathbb{Z}_{h}^{k}\brack m}\right|}{\left|{\mathbb{Z}_{h}^{n}\brack m}\right|}
=\prod_{i=1}^tp_i^{(s_i-1)(k-m)(n-k)}{n-m\brack k-m}_{p_i},$$
as desired.  $\qed$

Let $A$ be an $m$-subspace of $\mathbb{Z}_{h}^n$.  Let
$$A^\perp=\{y\in\mathbb{Z}_{h}^n: yx^t=0\;\hbox{for all}\;x\in A\}.$$
Then $A^\perp$ is a linear subset of $\mathbb{Z}_{h}^n$.  By Lemma~\ref{lemma2.2}, $A$ has a matrix representation
$A=(0_{m,n-m},I_m)T$, where $T\in G\!L_n(\mathbb{Z}_h)$. It is easy to prove that $A^\perp$  has a matrix
representation
\begin{equation}\label{equa20}
A^\perp=(I_{n-m},0_{n-m,m})(T^{-1})^{t}\; \hbox{if}\; A=(0_{m,n-m},I_m)T.
\end{equation}
Hence $A^\perp$ is an $(n-m)$-subspace of $\mathbb{Z}_{h}^n$.
The subspace $A^\perp$ is called the {\it dual subspace} of $A$. Note that
$\dim(A)+\dim(A^\perp)=n$ and $(A^{\perp})^\perp=A.$
If $A_1$ and $A_2$ are two subspaces of $\mathbb{Z}_{h}^n$, then $A_1\subseteq A_2$ if and only if $A_2^\perp\subseteq A_1^\perp.$

\begin{lemma}\label{lemma3.6}
Let $m<n$ and let $A$ and $B$ be two $m$-subspaces of $\mathbb{Z}_{h}^n$. Then
$$m-\dim(A\cap B)=n-m-\dim(A^\perp\cap B^\perp).$$
\end{lemma}
\proof Without loss of generality, by Lemma~\ref{lemma3.2}, we may assume that $A=(0_{m,n-m},I_m)$ and $B=(D,I_m)$,
where $D$ is as in (\ref{equa16}).
 By Theorem~\ref{lemma3.3}, $\dim(A\cap B)=m-r$.

 Note that  $A^\perp$ and $B^\perp$ are $(n-m)$-subspaces of $\mathbb{Z}_{h}^n$. By (\ref{equa20}), $A^\perp$ has
a matrix representation $A^\perp=(I_{n-m},0_{n-m,m})$. Let $B^\perp=(B_1,B_2)$ be a matrix representation
of $B^\perp$, where $B_1\in\mathbb{Z}_{h}^{(n-m)\times(n-m)}$ and $B_2\in\mathbb{Z}_{h}^{(n-m)\times m}$. Then $B_1D^t+B_2=0$,
which implies that $\rho(B_2)=\rho(-B_1D^t)\leq \rho(D)=r$ by (\ref{equa12}).
By (\ref{equa11}) and (\ref{equa17}), we have
\begin{eqnarray*}
n-m-\dim(A^\perp\cap B^\perp)&=&\rho\left(\begin{array}{c}
A^\perp\\ B^\perp
\end{array}\right)-(n-m)\\
&=&\rho\left(\begin{array}{cc}
I_{n-m} & 0_{n-m,m}\\ B_1 & B_2
\end{array}\right)-(n-m)\\
&=&\rho(B_2)\leq r=m-\dim(A\cap B).
\end{eqnarray*}
On the other hand,  we have similarly
$$m-\dim(A\cap B)=m-\dim((A^\perp)^\perp\cap(B^\perp)^\perp)
\leq n-m-\dim(A^\perp\cap B^\perp).$$
Therefore, we obtain $m-\dim(A\cap B)=n-m-\dim(A^\perp\cap B^\perp).$ $\qed$

\section{Proof of Theorem~\ref{thm1.1}}
In this section, we discuss the generalized Grassmann graph
$G_r(m,n,\mathbb{Z}_{h})$ over $\mathbb{Z}_{h}$ and prove Theorem~\ref{thm1.1}.

For convenience, if $X$ is a vertex of $G_r(m,n,\mathbb{Z}_{h})$, then the matrix representation of $X$ is also denoted by $X$.
By Lemmas~\ref{lemma2.10},~\ref{lemma2.11} and matrix representations of subspaces,  for each $i\in[t]$, $\pi_i$ induces a surjective  map by
\begin{equation}\label{equa21}
\begin{array}{rcl}
\pi_i: V(G_r(m,n,\mathbb{Z}_{h})) &\rightarrow & V(G_r(m,n,\mathbb{Z}_{p_i^{s_i}})),\\
X &\mapsto& \pi_i(X);
\end{array}
\end{equation}
and $\theta_i$ induces a surjective  map   by
\begin{equation}\label{equa22}
\begin{array}{rcl}
\theta_i: V(G_r(m,n,\mathbb{Z}_{h})) &\rightarrow & V(G_r(m,n,\mathbb{Z}_{h_i})),\\
X  &\mapsto&  \theta_i(X),
\end{array}
\end{equation}
where $h_i=h/p_i^{s_i}$.
Similarly, Lemmas~\ref{lemma2.10},~\ref{lemma2.11}, (\ref{equa6}) and (\ref{equa19}), $\pi$ induces a bijective  map by
\begin{equation}\label{equa23}
\begin{array}{rcl}
\pi: V(G_r(m,n,\mathbb{Z}_{h})) &\rightarrow& V(G_r(m,n,\mathbb{Z}_{p_1^{s_1}}))\times\cdots\times V(G_r(m,n,\mathbb{Z}_{p_t^{s_t}})),\\
X &\mapsto& (\pi_1(X),\ldots,\pi_t(X));
\end{array}
\end{equation}
and for each $i\in[t]$, $(\pi_i,\theta_i)$ induces a bijective  map  by
\begin{equation}\label{equa24}
\begin{array}{rcl}
(\pi_i,\theta_i): V(G_r(m,n,\mathbb{Z}_{h})) &\rightarrow& V(G_r(m,n,\mathbb{Z}_{p_i^{s_i}}))\times V(G_r(m,n,\mathbb{Z}_{h_i})),\\
X  &\mapsto&  (\pi_i(X),\theta_i(X)).
\end{array}
\end{equation}

\begin{thm}\label{lemma4.1}
Let $n\leq 2m$ and $2\leq r\leq \max\{m+1,n-m+1\}$. Then $G_r(m,n,\mathbb{Z}_{h})$ is a connected vertex-transitive graph, and
$G_r(m,n,\mathbb{Z}_{h})\cong G_r(n-m,n,\mathbb{Z}_{h})$.
\end{thm}
\proof If $X$ is a vertex of $G_r(m,n,\mathbb{Z}_{h})$, then there is a $T_X\in G\!L_n(\mathbb{Z}_h)$ such that
$X=(0_{m,n-m},I_m)T_X$ by Lemma~\ref{lemma2.2}. Let $\varphi(X)=XT_X^{-1}$. Then $\varphi$ is an automorphism
of $G_r(m,n,\mathbb{Z}_{h})$ and $\varphi(X)=(0_{m,n-m},I_m)$. Hence $G_r(m,n,\mathbb{Z}_{h})$ is vertex-transitive.

Let $A,B\in V(G_r(m,n,\mathbb{Z}_{h}))$ with $\dim(A\cap B)=m-\ell$. Without loss of generality, by Lemma~\ref{lemma3.2} and (\ref{equa17}),
we may assume that $A=(0_{m,n-m},I_m)$ and $B=(D,I_m)$,
where $$
D={\rm diag}\left(\prod_{i=1}^tp_i^{\alpha_{i1}},\prod_{i=1}^tp_i^{\alpha_{i2}},\ldots,\prod_{i=1}^tp_i^{\alpha_{i\ell}},0,\ldots,0\right)_{m,n-m},
$$
$1\leq \ell\leq \min\{m,n-m\},\rho(D)=\ell,$ and
 $0\leq\alpha_{i1}\leq\alpha_{i2}\leq\cdots\leq\alpha_{i\ell}\leq s_i$ for all $i\in[t]$.
For each $c\in[\ell]$, let $$A_c=(D_c,I_m)\;\hbox{with}\;D_c={\rm diag}\,\left(\prod_{i=1}^tp_i^{\alpha_{i1}},\prod_{i=1}^tp_i^{\alpha_{i2}},\ldots,
\prod_{i=1}^tp_i^{\alpha_{ic}},0,\ldots,0\right)_{m,n-m}.$$
Then $A_{\ell}=B$. Write $A_0=A$. Since $\rho\left(\begin{array}{c}A_{c-1}\\ A_c\end{array}\right)=m+1$ for each $c\in[\ell]$.
By (\ref{equa17}) and the definition of $G_r(m,n,\mathbb{Z}_{h})$, $A_{c-1}\sim A_c$ for each $c\in[\ell]$.
It follows that $G_r(m,n,\mathbb{Z}_{h})$ is a connected graph.

By Theorem~\ref{lemma3.5},  we have $|V(G_r(m,n,\mathbb{Z}_{h}))|=|V(G_r(n-m,n,\mathbb{Z}_{h}))|.$
By Lemma~\ref{lemma3.6}, the map $X\mapsto X^\perp$ is an isomorphism
from $G_r(m,n,\mathbb{Z}_{h})$ to $G_r(n-m,n,\mathbb{Z}_{h})$. Therefore, we obtain $G_r(m,n,\mathbb{Z}_{h})\cong G_r(n-m,n,\mathbb{Z}_{h})$. $\qed$

By Theorem~\ref{lemma4.1}, we may assume that $n\geq2m$ in our discussion on $G_r(m,n,\mathbb{Z}_{h})$.

\begin{lemma}\label{lemma4.2}{\rm (See Lemma~4.7 and Theorem~4.9 in \cite{Huang3}.)}
Let $n\geq 2m$. Then the clique number of $G_r(m,n,\mathbb{Z}_{p^s})$ is
$p^{(s-1)(n-m)(r-1)}{n-m+r-1\brack r-1}_{p}.$
Moreover, ${\cal F}$ is a maximum clique of $G_r(m,n,\mathbb{Z}_{p^s})$ if and only if
either {\rm(a)} ${\cal F}$ consists of all $m$-subspaces of $\mathbb{Z}_{p^s}^n$ which contain a fixed $(m-r+1)$-subspace
of $\mathbb{Z}_{p^s}^n$, or {\rm(b)} $n=2m$ and ${\cal F}$ is the set of all $m$-subspaces of $\mathbb{Z}_{p^s}^n$ contained in a fixed
$(m+r-1)$-subspace of $\mathbb{Z}_{p^s}^n$.
\end{lemma}

\begin{lemma}\label{lemma4.3}
Let $n\geq 2m$, and $\pi_i$ be as in (\ref{equa21}) for each $i\in[t]$.
Suppose that ${\cal F}$ is a maximum clique of $G_r(m,n,\mathbb{Z}_{h})$.  Then the following hold:
\begin{itemize}
\item[\rm(i)]
The clique number of $G_r(m,n,\mathbb{Z}_{h})$  is
$\prod_{i=1}^tp_i^{(s_i-1)(n-m)(r-1)}{n-m+r-1\brack r-1}_{p_i}.$

\item[\rm(ii)]
$\pi_i({\cal F})$ is a maximum clique of $G_r(m,n,\mathbb{Z}_{p_i^{s_i}})$ for all $i\in[t]$.
\end{itemize}
\end{lemma}
\proof  Let ${\cal F}$ be a maximum clique of $G_r(m,n,\mathbb{Z}_{h})$, and for each $i\in[t]$,
$$\pi_i({\cal F}):=\{\pi_i(F) : F\in{\cal F}\}=\{\pi_i(A_{i1}),\pi_i(A_{i2}),\ldots,\pi_i(A_{ik_i})\}.$$
By Lemma~\ref{lemma2.10}, (\ref{equa17-2}) and (\ref{equa21}), $\pi_i({\cal F})$ is a clique of $G_r(m,n,\mathbb{Z}_{p_i^{s_i}})$. 
For each $i\in[t]$, by Lemma~\ref{lemma4.2}, we have
\begin{equation}\label{equa26}
k_i\leq\omega(G_r(m,n,\mathbb{Z}_{p_i^{s_i}}))=p_i^{(s_i-1)(n-m)(r-1)}{n-m+r-1\brack r-1}_{p_i}.
\end{equation}
Moreover, ${\cal F}$ has a partition into $k_i$ cliques: ${\cal F}=\bigcup_{j=1}^{k_i}{\cal C}_j$, where ${\cal C}_j$ is a clique
with $\pi_i({\cal C}_j)=\{\pi_i(A_{ij})\}$ for each $j\in[k_i]=\{1,2,\ldots,k_i\}$, and ${\cal C}_u\cap{\cal C}_v=\emptyset$ for all $u\not=v$.
It follows that $\omega(G_r(m,n,\mathbb{Z}_{h}))=\sum_{j=1}^{k_i}|{\cal C}_j|$.

(i). By Lemma~\ref{lemma4.2}, the result is trivial for $t=1$. Suppose
that $t\geq2$ and the result is true for $t-1$.
Let $n_j=|{\cal C}_j|$ for each $j\in[k_i]$. Then there exists $\{B_{1j},B_{2j},\ldots,B_{n_jj}\}\subseteq J_{(0,\ldots,0,s_i,0,\ldots,0)}^{m\times n}$
such that ${\cal C}_j=\{A_{ij}+B_{1j},A_{ij}+B_{2j},\ldots,A_{ij}+B_{n_jj}\}$. By (\ref{equa24}),
$|\theta_i({\cal C}_j)|=|(\pi_i,\theta_i)({\cal C}_j)|=|{\cal C}_j|=n_j$.  By Lemma~\ref{lemma2.11}, (\ref{equa17-2}) and (\ref{equa22}), $\theta_i({\cal C}_j)$ is a clique of $G_r(m,n,\mathbb{Z}_{h_i})$, where $h_i=h/p_i^{s_i}$. By induction, for all $j\in[k_i]$, we have $$n_j\leq\omega(G_r(m,n,\mathbb{Z}_{h_i}))=\prod_{u\in[t]\setminus\{i\}}p_u^{(s_u-1)(n-m)(r-1)}{n-m+r-1\brack r-1}_{p_u}.$$
Therefore, by (\ref{equa26}), we have
\begin{equation}\label{equa27}
\omega(G_r(m,n,\mathbb{Z}_{h}))=\sum_{j=1}^{k_i}n_j\leq  \prod_{i=1}^tp_i^{(s_i-1)(n-m)(r-1)}{n-m+r-1\brack r-1}_{p_i}.
\end{equation}

On the other hand, $${\cal C}=\left\{\left(\begin{array}{cc}
X & 0\\
0 & I_{m-r+1}
\end{array}\right): X\in V(G_r(r-1,n-m+r-1,\mathbb{Z}_{h}))\right\}$$
is a clique of $G_r(m,n,\mathbb{Z}_{h})$ and $|{\cal C}|=\prod_{i=1}^tp_i^{(s_i-1)(n-m)(r-1)}{n-m+r-1\brack r-1}_{p_i}$.
It follows that $\omega(G_r(m,n,\mathbb{Z}_{h}))\geq\prod_{i=1}^tp_i^{(s_i-1)(n-m)(r-1)}{n-m+r-1\brack r-1}_{p_i}$. Therefore, the clique number of $G_r(m,n,\mathbb{Z}_{h})$  is $\prod_{i=1}^tp_i^{(s_i-1)(n-m)(r-1)}{n-m+r-1\brack r-1}_{p_i}$.

(ii). By (\ref{equa26}) and (\ref{equa27}), $k_i=p_i^{(s_i-1)(n-m)(r-1)}{n-m+r-1\brack r-1}_{p_i}$  for all $i\in[t]$. By Lemma~\ref{lemma4.2} again, $\pi_i(\cal F)$ is a maximum clique of $G_r(m,n,\mathbb{Z}_{p_i^{s_i}})$ for all $i\in[t]$.
 $\qed$

Let $4\leq2r\leq2m=n$ and ${\cal I}\subseteq[t]$. Then $2\leq2(m-r+1)\leq2(m-1)<2m=n$. Recall that
\begin{eqnarray*}
&&{\cal F}_{(\alpha_1,\alpha_2,\ldots,\alpha_t)}^{(m-r+1,m,n,{\cal I})}\\
&=&\left\{\left(\begin{array}{cc}
X & 0_{r-1,m-r+1}\\
Y & x_{{\cal I}}I_{m-r+1}
\end{array}\right)\in{\mathbb{Z}_h^n\brack m} :
X\in\mathbb{Z}_h^{(r-1)\times(m+r-1)},Y\in J_{(\alpha_1,\alpha_2,\ldots,\alpha_t)}^{(m-r+1)\times(m+r-1)} \right\},
\end{eqnarray*}
where $\pi_i(x_{\cal I})=1$ and $\alpha_i=s_i$ if $i\in {\cal I}$;
 and $\pi_i(x_{\cal I})=0$ and $\alpha_i=0$ if $i\in[t]\setminus {\cal I}$.

\begin{lemma}\label{lemma4.4}
Let $4\leq2r\leq2m=n$ and ${\cal I}\subseteq[t]$. Then ${\cal F}_{(\alpha_1,\alpha_2,\ldots,\alpha_t)}^{(m-r+1,m,n,{\cal I})}$
is a maximum clique of $G_r(m,n,\mathbb{Z}_{h})$.
\end{lemma}
\proof By (\ref{equa23}),  we have
$$\left|{\cal F}_{(\alpha_1,\alpha_2,\ldots,\alpha_t)}^{(m-r+1,m,n,{\cal I})}\right|=\left|\pi\left({\cal F}_{(\alpha_1,\alpha_2,\ldots,\alpha_t)}^{(m-r+1,m,n,{\cal I})}\right)\right|
=\prod_{i=1}^{t}\left|\pi_i\left({\cal F}_{(\alpha_1,\alpha_2,\ldots,\alpha_t)}^{(m-r+1,m,n,{\cal I})}\right)\right|.$$
For any
$$U=\left(\begin{array}{cc}
X & 0_{r-1,m-r+1}\\
Y & x_{{\cal I}}I_{m-r+1}
\end{array}\right)\in{\cal F}_{(\alpha_1,\alpha_2,\ldots,\alpha_t)}^{(m-r+1,m,n,{\cal I})},$$
 since $U\in{\mathbb{Z}_h^n\brack m}$, $X$ is an $(r-1)$-subspace of $\mathbb{Z}_h^{m+r-1}$.
By Theorem~\ref{lemma3.5} (i), one obtains $$\left|{\cal F}_{(\alpha_1,\alpha_2,\ldots,\alpha_t)}^{(m-r+1,m,n,{\cal I})}\right|
\geq\left|{\mathbb{Z}_h^{m+r-1}\brack r-1}\right|=\prod_{i=1}^tp_i^{(s_i-1)m(r-1)}{m+r-1\brack r-1}_{p_i}.$$
If $i\in {\cal I}$, from Lemma~\ref{lemma2.10}, we deduce that $\pi_i\left({\cal F}_{(\alpha_1,\alpha_2,\ldots,\alpha_t)}^{(m-r+1,m,n,{\cal I})}\right)\subseteq{\cal A}_i$, where ${\cal A}_i$ is the set of all $m$-subspaces of $\mathbb{Z}_{p_i^{s_i}}^n$
containing the fixed $(m-r+1)$-subspace $(0_{m-r+1,m+r-1},I_{m-r+1})$.
If $i\in[t]\setminus {\cal I}$, from Lemma~\ref{lemma2.10}, we deduce that $\pi_i\left({\cal F}_{(\alpha_1,\alpha_2,\ldots,\alpha_t)}^{(m-r+1,m,n,{\cal I})}\right)\subseteq{\cal B}_i$,
where ${\cal B}_i$ is the set of all $m$-subspaces of $\mathbb{Z}_{p_i^{s_i}}^n$
contained in the fixed $(m+r-1)$-subspace $(I_{m+r-1},0_{m+r-1,m-r+1})$. Therefore, by Theorem~\ref{lemma3.5} (i) again, we find
$$\prod_{i=1}^{t}\left|\pi_i\left({\cal F}_{(\alpha_1,\alpha_2,\ldots,\alpha_t)}^{(m-r+1,m,n,{\cal I})}\right)\right|
\leq\prod_{i\in{\cal I}}|{\cal A}_i|\prod_{i\in[t]\setminus{\cal I}}|{\cal B}_i|
=\prod_{i=1}^tp_i^{(s_i-1)m(r-1)}{m+r-1\brack r-1}_{p_i}.$$
Hence, we obtain \begin{equation}\label{equa28}
\left|{\cal F}_{(\alpha_1,\alpha_2,\ldots,\alpha_t)}^{(m-r+1,m,n,{\cal I})}\right|=\prod_{i=1}^tp_i^{(s_i-1)m(r-1)}{m+r-1\brack r-1}_{p_i}.
\end{equation}

Let $U$ and $W$ be two distinct subspaces in ${\cal F}_{(\alpha_1,\alpha_2,\ldots,\alpha_t)}^{(m-r+1,m,n,{\cal I})}$.
If $i\in {\cal I}$, then $\pi_i(U)\cap \pi_i(W)\supseteq(0_{m-r+1,m+r-1},I_{m-r+1})$, which implies that $\dim(\pi_i(U)\cap \pi_i(W))\geq m-r+1$.
If $i\in[t]\setminus {\cal I}$, then $\pi_i(U),\pi_i(W)\subseteq(I_{m+r-1},0_{m+r-1,m-r+1})$, by (\ref{equa17}) and $\dim(\pi_i(U)\vee\pi_i(W))\leq m+r-1$,
we have $\dim(\pi_i(U)\cap \pi_i(W))\geq m-r+1$.
Therefore, we have $\dim(\pi_i(U)\cap \pi_i(W))\geq m-r+1$ for all $i\in[t]$. By (\ref{equa17-2}), we have $\dim(U\cap W)\geq m-r+1$,
and therefore ${\cal F}_{(\alpha_1,\alpha_2,\ldots,\alpha_t)}^{(m-r+1,m,n,{\cal I})}$ is a clique of $G_r(m,n,\mathbb{Z}_{h})$.
By Lemma~\ref{lemma4.3} and (\ref{equa28}),
${\cal F}_{(\alpha_1,\alpha_2,\ldots,\alpha_t)}^{(m-r+1,m,n,{\cal I})}$ is a maximum clique of $G_r(m,n,\mathbb{Z}_{h})$. $\qed$

\begin{thm}\label{lemma4.5}
Let $n\geq 2m$. Then ${\cal F}$ is a maximum clique of $G_r(m,n,\mathbb{Z}_{h})$ if and only if
either {\rm(a)} ${\cal F}$ consists of all $m$-subspaces of $\mathbb{Z}_{h}^n$ which contain a fixed $(m-r+1)$-subspace
of $\mathbb{Z}_{h}^n$,  {\rm(b)} $n=2m$ and ${\cal F}$ is the set of all $m$-subspaces of $\mathbb{Z}_{h}^n$ contained in a fixed
$(m+r-1)$-subspace of $\mathbb{Z}_{h}^n$, or {\rm(c)} $n=2m$ and there exists some $T\in G\!L_n(\mathbb{Z}_h)$ such that 
${\cal F}=\left\{FT : F\in{\cal F}_{(\alpha_1,\alpha_2,\ldots,\alpha_t)}^{(m-r+1,m,n,{\cal I})}\right\}$
with ${\cal I}\not=\emptyset,[t]$.
\end{thm}
\proof By the definition of $G_r(m,n,\mathbb{Z}_{h})$, we have $2\leq r\leq m+1$. When $r=m+1$ or $m=1$, $G_r(m,n,\mathbb{Z}_{h})$ is a
clique, which implies that this theorem holds by Theorem~\ref{lemma3.5}. From now on, we assume that $4\leq2r\leq2m\leq n$.

Suppose that either {\rm(a)} ${\cal F}$ consists of all $m$-subspaces of $\mathbb{Z}_{h}^n$ which contain a fixed $(m-r+1)$-subspace
of $\mathbb{Z}_{h}^n$,  {\rm(b)} $n=2m$ and ${\cal F}$ is the set of all $m$-subspaces of $\mathbb{Z}_{h}^n$ contained in a fixed
$(m+r-1)$-subspace of $\mathbb{Z}_{h}^n$, or {\rm(c)} $n=2m$ and  there exists some $T\in G\!L_n(\mathbb{Z}_h)$ such that ${\cal F}={\cal F}_{(\alpha_1,\alpha_2,\ldots,\alpha_t)}^{(m-r+1,m,n,{\cal I})}T$
with ${\cal I}\not=\emptyset,[t]$. By the definition of $G_r(m,n,\mathbb{Z}_{h})$, Theorem~\ref{lemma3.5}, Lemma~\ref{lemma4.3} and Lemma~\ref{lemma4.4},  ${\cal F}$ is a maximum clique of $G_r(m,n,\mathbb{Z}_{h})$.

Let ${\cal F}$ be a maximum clique of $G_r(m,n,\mathbb{Z}_{h})$.  We will prove that either {\rm(a)} ${\cal F}$ consists of all $m$-subspaces of $\mathbb{Z}_{h}^n$ which contain a fixed $(m-r+1)$-subspace of $\mathbb{Z}_{h}^n$,  {\rm(b)} $n=2m$ and ${\cal F}$ is the set of all $m$-subspaces of $\mathbb{Z}_{h}^n$ contained in a fixed
$(m+r-1)$-subspace of $\mathbb{Z}_{h}^n$, or {\rm(c)} $n=2m$ and  there exists some $T\in G\!L_n(\mathbb{Z}_h)$ such that ${\cal F}={\cal F}_{(\alpha_1,\alpha_2,\ldots,\alpha_t)}^{(m-r+1,m,n,{\cal I})}T$
with ${\cal I}\not=\emptyset,[t]$. There are the following two cases to be considered.

{\bf Case}~1: $n>2m$. For each $i\in [t]$, by Lemmas~\ref{lemma4.2} and~\ref{lemma4.3},  $\pi_i({\cal F})$ consists of all $m$-subspaces of $\mathbb{Z}_{p_i^{s_i}}^n$ which contain a fixed $(m-r+1)$-subspace $P_i$ of $\mathbb{Z}_{p_i^{s_i}}^n$. By Lemma~\ref{lemma2.2}, there is a $T_i\in G\!L_n(\mathbb{Z}_{p_i^{s_i}})$ such that
$P_i=(0_{m-r+1,n-m+r-1},I_{m-r+1})T_i$. So,  we have
$$\pi_i({\cal F})=\left\{\left(\begin{array}{cc}
X_i & 0_{r-1,m-r+1}\\
0_{m-r+1,n-m+r-1} & I_{m-r+1}
\end{array}\right)T_i\in{\mathbb{Z}_{p_i^{s_i}}^{n}\brack m} : X_i\in{\mathbb{Z}_{p_i^{s_i}}^{n-m+r-1}\brack r-1}\right\}.$$
By Lemma~\ref{lemma2.6}, there is the unique $T\in G\!L_n(\mathbb{Z}_{h})$ such that $\pi(T)=\pi^\ast(T)=(T_1,\ldots,T_t)$.
For a given $(X_1,\ldots,X_t)\in{\mathbb{Z}_{p_1^{s_1}}^{n-m+r-1}\brack r-1}\times\cdots\times{\mathbb{Z}_{p_t^{s_t}}^{n-m+r-1}\brack r-1}$, from (\ref{equa23}), we deduce that
there is the unique $X\in{\mathbb{Z}_{h}^{n-m+r-1}\brack r-1}$ such that $\pi(X)=(X_1,\ldots,X_t)$. Therefore, for all $i\in[t]$, by (\ref{equa7}), we obtain
$$\pi_i({\cal F})=\left\{\pi_i\left(\left(\begin{array}{cc}
X & 0_{r-1,m-r+1}\\
0_{m-r+1,n-m+r-1} & I_{m-r+1}
\end{array}\right)T\right)\in{\mathbb{Z}_{p_i^{s_i}}^{n}\brack m} : X\in{\mathbb{Z}_{h}^{n-m+r-1}\brack r-1}\right\}.$$
It follows that
 $$\pi({\cal F})=\left\{\pi\left(\left(\begin{array}{cc}
X & 0_{r-1,m-r+1}\\
0_{m-r+1,n-m+r-1} & I_{m-r+1}
\end{array}\right)T\right) : X\in{\mathbb{Z}_{h}^{n-m+r-1}\brack r-1}\right\}.$$
From (\ref{equa23}) again, we deduce that
 $${\cal F}=\left\{\left(\begin{array}{cc}
X & 0_{r-1,m-r+1}\\
0_{m-r+1,n-m+r-1} & I_{m-r+1}
\end{array}\right)T\in{\mathbb{Z}_{h}^{n}\brack m} : X\in{\mathbb{Z}_{h}^{n-m+r-1}\brack r-1}\right\},$$
which implies that ${\cal F}$ consists of all $m$-subspaces of $\mathbb{Z}_{h}^n$ which contain a fixed $(m-r+1)$-subspace
of $\mathbb{Z}_{h}^n$.

{\bf Case}~2: $n=2m$. By Lemmas~\ref{lemma4.2} and~\ref{lemma4.3}, there exists some subset ${\cal I}$ of $[t]$ such that,
$\pi_i({\cal F})$ consists of all $m$-subspaces of $\mathbb{Z}_{p_i^{s_i}}^n$ which contain a fixed $(m-r+1)$-subspace of $\mathbb{Z}_{p_i^{s_i}}^n$ if $i\in {\cal I}$,   and $\pi_i({\cal F})$ is the set of all $m$-subspaces of $\mathbb{Z}_{p_i^{s_i}}^n$ contained in a fixed
$(m+r-1)$-subspace of $\mathbb{Z}_{p_i^{s_i}}^n$ if $i\in[t]\setminus {\cal I}$.

{\it Case}~2.1: ${\cal I}=[t]$. Similar to the proof of Case~1, we may obtain that
${\cal F}$ consists of all $m$-subspaces of $\mathbb{Z}_{h}^n$ which contain a fixed $(m-r+1)$-subspace
of $\mathbb{Z}_{h}^n$.

{\it Case}~2.2: ${\cal I}=\emptyset$. For each $i\in[t]$, $\pi_i({\cal F})$ is the set of all $m$-subspaces of $\mathbb{Z}_{p_i^{s_i}}^n$ contained in a fixed $(m+r-1)$-subspace $Q_i$ of $\mathbb{Z}_{p_i^{s_i}}^n$. By Lemma~\ref{lemma2.2}, there is a $T_i\in G\!L_n(\mathbb{Z}_{p_i^{s_i}})$ such that $Q_i=(I_{m+r-1},0_{m+r-1,m-r+1})T_i$. So,  we have
$$\pi_i({\cal F})=\left\{(X_i, 0_{m,m-r+1})T_i\in{\mathbb{Z}_{p_i^{s_i}}^{n}\brack m} : X_i\in{\mathbb{Z}_{p_i^{s_i}}^{m+r-1}\brack m}\right\}.$$
Let $T\in G\!L_n(\mathbb{Z}_{h})$ such that $\pi(T)=(T_1,\ldots,T_t)$.  
For all $i\in[t]$, by (\ref{equa7}) and  (\ref{equa23}), we have
$$\pi_i({\cal F})=\left\{\pi_i((X,0_{m,m-r+1})T)\in{\mathbb{Z}_{p_i^{s_i}}^{n}\brack m} : X\in{\mathbb{Z}_{h}^{m+r-1}\brack m}\right\}.$$
It follows that
$$\pi({\cal F})=\left\{\pi((X,0_{m,m-r+1})T) : X\in{\mathbb{Z}_{h}^{m+r-1}\brack m}\right\}.$$
Similarly, one obtains
 $${\cal F}=\left\{(X,0_{m,m-r+1})T\in{\mathbb{Z}_{h}^{n}\brack m} : X\in{\mathbb{Z}_{h}^{m+r-1}\brack m}\right\}.$$
Thus, ${\cal F}$ is the set of all $m$-subspaces of $\mathbb{Z}_{h}^n$ contained in a fixed
$(m+r-1)$-subspace $(I_{m+r-1},0_{m+r-1,m-r+1})T$ of $\mathbb{Z}_{h}^n$.

{\it Case}~2.3: ${\cal I}\not=\emptyset,[t]$. If $i\in {\cal I}$,
$\pi_i({\cal F})$ consists of all $m$-subspaces of $\mathbb{Z}_{p_i^{s_i}}^n$ which contain a fixed $(m-r+1)$-subspace $P_i$ of $\mathbb{Z}_{p_i^{s_i}}^n$.
By Lemma~\ref{lemma2.2}, there is a $T_i\in G\!L_n(\mathbb{Z}_{p_i^{s_i}})$ such that $P_i=(0_{m-r+1,m+r-1},I_{m-r+1})T_i$. So,  we have
$$\pi_i({\cal F})=\left\{\left(\begin{array}{cc}
X_i & 0_{r-1,m-r+1}\\
0_{m-r+1,m+r-1} & I_{m-r+1}
\end{array}\right)T_i\in{\mathbb{Z}_{p_i^{s_i}}^{n}\brack m} : X_i\in{\mathbb{Z}_{p_i^{s_i}}^{m+r-1}\brack r-1}\right\}.$$

If $i\in[t]\setminus {\cal I}$, $\pi_i({\cal F})$ is the set of all $m$-subspaces of $\mathbb{Z}_{p_i^{s_i}}^n$ contained in a fixed $(m+r-1)$-subspace $Q_i$ of $\mathbb{Z}_{p_i^{s_i}}^n$. By Lemma~\ref{lemma2.2},  there is a $T_i\in G\!L_n(\mathbb{Z}_{p_i^{s_i}})$ such that $Q_i=(I_{m+r-1},0_{m+r-1,m-r+1})T_i$. So,  we have
$$\pi_i({\cal F})=\left\{\left(\begin{array}{cc}
X_i & 0_{r-1,m-r+1}\\
Y_i & 0_{m-r+1,m-r+1}
\end{array}\right)T_i\in{\mathbb{Z}_{p_i^{s_i}}^{n}\brack m} :
\left(\begin{array}{c}
X_i\\ Y_i
\end{array}\right)\in{\mathbb{Z}_{p_i^{s_i}}^{m+r-1}\brack m}\right\}.$$

Let $T\in G\!L_n(\mathbb{Z}_{h})$ satisfy $\pi(T)=(T_1,\ldots,T_t)$. From (\ref{equa7}) and (\ref{equa23}), we deduce that for all $i\in{\cal I}$, 
$$
\pi_i({\cal F})\subseteq
\left\{\pi_i\left(\left(\begin{array}{cc}
X & 0_{r-1,m-r+1}\\
0_{m-r+1,m+r-1} & I_{m-r+1}
\end{array}\right)T\right)\in{\mathbb{Z}_{p_i^{s_i}}^{n}\brack m} : X\in{\mathbb{Z}_{h}^{m+r-1}\brack r-1}\right\};  
$$
and for all $i\in[t]\setminus{\cal I}$, 
$$\pi_i({\cal F})\subseteq\left\{\pi_i\left(\left(\begin{array}{cc}
X & 0_{r-1,m-r+1}\\
Y & 0_{m-r+1,m-r+1}
\end{array}\right)T\right)\in{\mathbb{Z}_{p_i^{s_i}}^{n}\brack m} :
\left(\begin{array}{c}
X\\ Y
\end{array}\right)\in{\mathbb{Z}_{h}^{m+r-1}\brack m}\right\}.$$

Note that $\pi_i(J_{(\alpha_1,\alpha_2,\ldots,\alpha_t)})=(0)$ if $\alpha_i=s_i$, and $\pi_i(J_{(\alpha_1,\alpha_2,\ldots,\alpha_t)})=\mathbb{Z}_{p_i^{s_i}}$ if $\alpha_i=0$. 
Let  $\pi_i(x_{\cal I})=1$ and $\alpha_i=s_i$ if $i\in {\cal I}$;
 and $\pi_i(x_{\cal I})=0$ and $\alpha_i=0$ if $i\in[t]\setminus {\cal I}$. For all $i\in[t]$, we have
\begin{eqnarray*}
\pi_i({\cal F})&\subseteq&\left\{\pi_i\left(\left(\begin{array}{cc}
X & 0_{r-1,m-r+1}\\
Y & x_{\cal I}I_{m-r+1}
\end{array}\right)T\right)\in{\mathbb{Z}_{p_i^{s_i}}^n \brack m} :\right.\\
&&\left.X\in\mathbb{Z}_h^{(r-1)\times(m+r-1)},Y\in J_{(\alpha_1,\alpha_2,\ldots,\alpha_t)}^{(m-r+1)\times(m+r-1)}\right\}.
\end{eqnarray*}
It follows that $\pi({\cal F})\subseteq\left\{\pi(FT): F\in{\cal F}_{(\alpha_1,\alpha_2,\ldots,\alpha_t)}^{(m-r+1,m,n,{\cal I})}\right\}$. So,
${\cal F}\subseteq\left\{FT: F\in{\cal F}_{(\alpha_1,\alpha_2,\ldots,\alpha_t)}^{(m-r+1,m,n,{\cal I})}\right\}$. 
Since ${\cal F}$ is a maximum clique  and $\left\{FT: F\in{\cal F}_{(\alpha_1,\alpha_2,\ldots,\alpha_t)}^{(m-r+1,m,n,{\cal I})}\right\}$ 
is a clique of $G_{r}(m,n,\mathbb{Z}_{h})$,  ${\cal F}=\left\{FT: F\in{\cal F}_{(\alpha_1,\alpha_2,\ldots,\alpha_t)}^{(m-r+1,m,n,{\cal I})}\right\}$.
$\qed$

\medskip\noindent{\bf Proof of Theorem~\ref{thm1.1}.}
Let $\lfloor n/2\rfloor\geq m\geq r\geq 0$ and ${\cal F}\subseteq{\mathbb{Z}_{h}^{n}\brack m}$ be an $r$-intersecting family. Without loss of
generality, we assume that $1\leq r\leq m-1$. By Theorem~\ref{lemma3.3}, ${\cal F}$ is a clique of the generalized Grassmann graph
$G_{m-r+1}(m,n,\mathbb{Z}_{h})$. By Lemma~\ref{lemma4.3}, $|{\cal F}|\leq\prod_{i=1}^tp_i^{(s_i-1)(n-m)(m-r)}{n-r\brack m-r}_{p_i}$
  and equality holds if and only if ${\cal F}$ is a
maximum clique of $G_{m-r+1}(m,n,\mathbb{Z}_{h})$. By Theorem~\ref{lemma4.5}, this theorem holds. $\qed$

\section*{Acknowledgment}
This research is supported by
 National Natural Science Foundation of China (11971146).

\end{document}